\documentclass[11pt,reqno]{amsart}
\usepackage{amsmath,amsfonts, amssymb, amsthm}
  \usepackage{paralist}
  \usepackage{graphics} 
  \usepackage{epsfig} 
\usepackage{graphicx}  \usepackage{epstopdf}
 \usepackage[colorlinks=true]{hyperref}
\hypersetup{urlcolor=blue, citecolor=red}

\newtheorem{theorem}{Theorem}[section]
\newtheorem{corollary}{Corollary}[section]

\newtheorem{lemma}[theorem]{Lemma}

\newtheorem{conjecture}{Conjecture}[section]

\theoremstyle{definition}

\newtheorem{remark}{Remark}[section]

\def\pmod #1{\ ({\rm{mod}}\ #1)}
\def\Z{\Bbb Z}
\def\N{\Bbb N}

\def\C{\Bbb C}

\def\l{\left}
\def\r{\right}
\def\bg{\bigg}
\def\({\bg(}
\def\){\bg)}
\def\t{\text}
\def\f{\frac}
\def\mo{{\rm{mod}\ }}
\def\pmod#1{\ (\mo\ #1)}

\def\ls{\leq}
\def\gs{\geq}

\def\al{\alpha}

\def\ve{\varepsilon}

\def\eq{\equiv}

\def\da{\delta}
\def\la{\lambda}

\def\Proof{\noindent{\it Proof}}

\begin{document}
\hbox{Preprint}
\medskip

\title[On determinants involving second-order recurrent sequences]
      {On determinants involving second-order recurrent sequences}
\author[Zhi-Wei Sun]{Zhi-Wei Sun}


\address{Department of Mathematics, Nanjing
University, Nanjing 210093, People's Republic of China}
\email{{\tt zwsun@nju.edu.cn}
\newline\indent
{\it Homepage}: {\tt http://maths.nju.edu.cn/\lower0.5ex\hbox{\~{}}zwsun}}

\keywords{Second-order recurrence, Lucas sequence, Fibonacci number, determinant, tree.
\newline \indent Supported by the Natural Science Foundation of China (grant no. 11971222).}

\begin{abstract}
Let $A$ and $B$ be complex numbers, and let $(w_n)_{n\gs0}$ be a sequence of complex numbers
with $w_{n+1}=Aw_n-Bw_{n-1}$ for all $n=1,2,3,\ldots$. When $w_0=0$ and $w_1=1$, the sequence $(w_n)_{n\gs0}$ is just the Lucas sequence $(u_n(A,B))_{n\gs0}$. In this paper, we evaluate the determinants
$$\det[w_{|j-k|}]_{1\ls j,k\ls n}\ \ \t{and}\ \ \det[w_{|j-k+1|}]_{1\ls j,k\ls n}.$$
In particular, we have
$$\det[u_{|j-k|}(A,B)]_{1\ls j,k\ls n}=(-1)^{n-1}u_{n-1}(2A,(B+1)^2).$$
When $B=-1$ and $2\mid n$, we also determine the characteristic polynomial of the matrix $[w_{j+k}]_{0\ls j,k\ls n-1}$.
\end{abstract}
\maketitle

\section{Introduction}
\setcounter{equation}{0}
 \setcounter{conjecture}{0}
 \setcounter{theorem}{0}
 \setcounter{proposition}{0}

 In 1934 R. Robinson proposed the evaluation of the determinant $\det[|j-k|]_{1\ls j,k\ls n}$ as a problem in Amer. Math. Monthly, later its solutions appeared in \cite{RS}. Namely, we have
\begin{equation}\label{j-k}\det[|j-k|]_{1\ls j,k\ls n}=(-1)^{n-1}(n-1)2^{n-2}.
\end{equation}

Let $n\in\Z^+=\{1,2,3,\ldots\}$, and let $T$ be any (undirected) tree with $n$ vertices $v_1,\ldots,v_n$. For $j,k=1,\ldots,n$, let $d(v_j,v_k)$ denote the distance between the vertices $v_j$ and $v_k$.
In 1971 R.L. Graham and H.O. Pollak \cite{GP} established the following celebrated formula:
\begin{equation}\label{tree}\det[d(v_j,v_k)]_{1\ls j,k\ls n}=(-1)^{n-1}(n-1)2^{n-2}.
\end{equation}
This is a further extension of \eqref{j-k} as a path with $n$ vertices is a tree.
Based on the idea in \cite{YY1}, in 2007 W. Yan and Y.-N. Yeh \cite[Corollary 2.3] {YY2} obtained the following $q$-analogue of \eqref{tree}
 for $n>1$:
\begin{equation}\label{q-tree}\det[[d(v_j,v_k)]_q]_{1\ls j,k\ls n}=(-1)^{n-1}(n-1)(1+q)^{n-2},
\end{equation}
where $[m]_q$ with $m\in\N$ denotes the $q$-analogue of $m$ given by
$$[m]_q:=\sum_{0\ls k<m}q^k=\begin{cases}(1-q^m)/(1-q)&\t{if}\ q\not=1,\\m&\t{if}\ q=1.\end{cases}$$
(Throughout this paper, we consider $0^0$ as $1$.)
Another result of Yan and Yeh \cite[Corollary 2.2]{YY2} states that
\begin{equation}\label{q}\det[q^{d(v_j,v_k)}]_{1\ls j,k\ls n}=(1-q^2)^{n-1},
\end{equation}

Let $R$ be a commutative ring with identity.
The Lucas sequences $(u_n(x,y))_{n\gs0}$ and $(v_n(x,y))_{n\gs0}$ over $R$ are defined as follows:
\begin{align*}&u_0(x,y)=0,\ u_1(x,y)=1,\ u_{n+1}(x,y)=xu_n(x,y)-yu_{n-1}(x,y)\ \t{for}\ n\in\Z^+;
\\&v_0(x,y)=2,\ v_1(x,y)=x,\ u_{n+1}(x,y)=xu_n(x,y)-yu_{n-1}(x,y)\ \t{for}\ n\in\Z^+.
\end{align*} It is well known that
$$u_n(x,y)=\f{x^n-y^n}{x-y}=\sum_{0\ls k<n}x^ky^{n-1-k}\ \ \t{and}\ \ v_n(x,y)=x^n+y^n$$
for all $n\in\N$. Note that for any $n\in\N$ we have $u_n(q+1,q)=[n]_q$, in particular $u_n(2,1)=n$.

Let $A,B\in\C$, where $\C$ is the field of complex numbers.
Let $(w_n)_{n\gs0}$ be a sequence of complex numbers satisfying the recurrence
\begin{equation}\label{w}w_{n+1}=Aw_n-B_{n-1}\ (n=1,2,3,\ldots).
\end{equation}
When $w_0=0$ and $w_1=1$, we have $w_n=u_n(A,B)$ for all $n\in\N$.
When $w_0=2$ and $w_1=A$, we have $w_n=v_n(A,B)$ for any $n\in\N$.
In this paper, we evaluate
$$\det[w_{|j-k|}]_{1\ls j,k\ls n}\ \ \t{and}\ \ \det[w_{|j-k+1|}]_{1\ls j,k\ls n},$$
which extends \eqref{j-k} in a new way.

\begin{theorem}\label{Th1.1} Let $A$ and $B$ be elements of a commutative ring $R$ with identity, and let $(w_n)_{n\gs0}$ be a sequence of elements of $R$
satisfying the recurrence \eqref{w}. For any $n\in\Z^+=\{1,2,3,\ldots\}$, we have
\begin{equation}\label{w'}\det[w_{|j-k|}]_{1\ls j,k\ls n}=w_0u_n(A',B')+((Bw_0)^2-(Aw_0-w_1)^2)u_{n-1}(A',B'),
\end{equation}
where
$$A'=(A^2-B^2+1)w_0-2Aw_1\ \t{and}\ B'=((Aw_0-(B+1)w_1)^2.$$
\end{theorem}

Taking $w_0=0$ and $w_1=1$ in Theorem \ref{Th1.1}
and noting that
\begin{equation}\label{xyz}u_n(xz,yz^2)=u_n(x,y)z^{n-1}\ \ (n=1,2,3,\ldots),
\end{equation}
we immediately obtain the following corollary.

\begin{corollary} Let $A$ and $B$ be elements of a commutative ring $R$ with identity. Then,
for any positive integer $n$, we have
\begin{equation}\label{u}\det[u_{|j-k|}(A,B)]_{1\ls j,k\ls n}=(-1)^{n-1}u_{n-1}(2A,(B+1)^2).
\end{equation}
\end{corollary}

Let $n>1$ be an integer. In view of \eqref{xyz},
$$u_{n-1}(2A,4)=2^{n-2}u_{n-1}(A,1).$$
So, \eqref{u} with $B=1$ gives the identity
\begin{equation}\det[u_{|j-k|}(A,1)]_{1\ls j,k\ls n}=(-1)^{n-1}2^{n-2}u_{n-1}(A,1).
\end{equation}
In the case $A=2$, this turns out to be the classical formula \eqref{j-k}.
Note also that \eqref{u} with $B=-1$ yields the identity
\begin{equation}\det[u_{|j-k|}(A,-1)]_{1\ls j,k\ls n}=(-1)^{n-1}(2A)^{n-2}.
\end{equation}
The identity \eqref{u} with $A=1$ and $B=-2$ gives the formula
\begin{equation}\det[u_{|j-k|}(1,-2)]_{1\ls j,k\ls n}=(-1)^{n-1}(n-1).
\end{equation}
In the case $B=q$ and $A=q+1$, with the aid of \eqref{xyz}, from the identity \eqref{u} we obtain the $q$-analogue of \eqref{j-k}:
\begin{equation}\det[[|j-k|]_q]_{1\ls j,k\ls n}=(-1)^{n-1}(n-1)(q+1)^{n-2}.
\end{equation}

One may wonder whether the identity \eqref{u} can be extended to trees.
The answer is negative. Let's consider a tree $T$ with vertices $v_1,v_2,v_3,v_4$
and edges $v_1v_2,v_2v_3,v_2v_4$. For any $A,B\in\C$, we clearly have
$$\det\l[u_{d(v_j,v_k)}(A,B)\r]_{1\ls j,k\ls 4}
=\begin{vmatrix}
0&1&A&A\\1&0&1&1\\A&1&0&A\\A&1&A&0
\end{vmatrix}=-3A^2$$
which is independent of $B$, while the right-hand side of \eqref{u} indeed depends on $B$.

It is easy to see that
$$u_n(1,1)=(-1)^{n-1}\l(\f n3\r)\quad\t{for all}\ n\in\N,$$
where $(-)$ denotes the Legendre symbol.
In view of this and \eqref{xyz}, the identity \eqref{u} with $B=\pm 2A-1$ yields the following corollary.

\begin{corollary} Let $A$ be any element of a commutative ring with identity. For any integer $n\gs2$ we have
\begin{equation}\begin{aligned}\det[u_{|j-k|}(A,2A-1)]_{1ls j,k\ls n}
&=\det[u_{|j-k|}(A,-2A-1)]_{1ls j,k\ls n}
\\&=\l(\f{1-n}3\r)(2A)^{n-2}.\end{aligned}
\end{equation}
\end{corollary}

Taking $w_0=2$ and $w_1=A$ in Theorem \ref{Th1.1} and making use of \eqref{xyz}, we get the following corollary.

\begin{corollary} \label{Cor-v} Let $A$ and $B$ be elements of a commutative ring $R$ with identity. For any integer $n>1$, we have
\begin{equation}\begin{aligned}\det[v_{|j-k|}(A,B)]_{1\ls j,k\ls n}=\ &2(1-B)^{n-1}u_n(2(1+B),A^2)
\\&\ +(4B^2-A^2)(1-B)^{n-2}u_{n-1}(2(1+B),A^2).
\end{aligned}\end{equation}
\end{corollary}

By Corollary \ref{Cor-v}, for any $A\in\C$, we have
\begin{equation}\det[v_{|j-k|}(A,1)]_{1\ls j,k\ls n}=0\quad \t{for all}\ n=3,4,5,\ldots.
\end{equation}
In the case $B=-1$, Corollary \ref{Cor-v} yields the following result.

\begin{corollary} For any $A\in\C$ and $n\in\{2,3,4,\ldots\}$, we have
\begin{equation}\det[v_{|j-k|}(A,-1)]_{1\ls j,k\ls n}
=(-1)^{\lfloor(n-1)/2\rfloor}(2A)^{n-2}\times\begin{cases}4A&\t{if}\ 2\nmid n,
\\4-A^2&\t{if}\ 2\mid n.\end{cases}
\end{equation}
\end{corollary}

Applying Corollary \ref{Cor-v} with $A=\pm2B$ and using the identity \eqref{xyz}, we obtain the following corollary.

\begin{corollary} Let $R$ be a commutative ring with identity. For any $B\in R$ and $n\in\Z^+$, we have
\begin{equation}\begin{aligned}\det[v_{|j-k|}(2B,B)]_{1\ls j,k\ls n}
&=\det[v_{|j-k|}(-2B,B)]_{1\ls j,k\ls n}
\\&=2^n(1-B)^{n-1}u_n(1+B,B^2).
\end{aligned}
\end{equation}
\end{corollary}

Recall that the Fibonacci numbers are those $F_n=u_n(1,-1)$ with $n\in\N$.
For any $n\in\N$, we clearly have
$$F_{2n+2}=F_{2n}+(F_{2n}+F_{2n-1})=2F_{2n}+(F_{2n}-F_{2n-2})=3F_{2n}-F_{2n-2}.$$
Thus $F_{2n}=u_n(3,1)$ for all $n\in\N$.

\begin{corollary} For any integer $n\gs2$, we have
\begin{equation}\label{u(3,1)}\det[v_{|j-k|}(2,2)]_{1\ls j,k\ls n}=(-2)^nF_{2n-4}.
\end{equation}
\end{corollary}
\Proof. In view of Corollary \ref{Cor-v} with $A=B=2$ and the identity \eqref{xyz},
\begin{align*}\det[v_{|j-k|}(2,2)]_{1\ls j,k\ls n}&=(-1)^{n-1}2(u_n(6,4)-6u_{n-1}(6,4))
\\&=(-1)^{n-1}2(-4u_{n-2}(6,4))
\\&=-8(-1)^{n-1}2^{n-3}u_{n-2}(3,1)=(-2)^nu_{n-2}(3,1).
\end{align*}
This implies \eqref{u(3,1)} since $u_m(3,1)=F_{2m}$ for all $m\in\N$. \qed

\begin{corollary} Let $R$ be a commutative ring with identity,
and let $A,B\in R$ and $\ve\in\{\pm1\}$.  Suppose that
$$w_{-1}=\ve,\ w_0=1,\ \t{and}\ w_{n+1}=Aw_n-Bw_{n-1}\ \t{for}\ n\in\N.$$
Then, for any $n\in\Z^+$, we have
\begin{equation}\label{sw}\det[w_{|j-k|}]_{1\ls j,k\ls n}=u_n(1-(A-\ve B)^2,B^2(1+B-\ve A)^2).
\end{equation}
\end{corollary}
\Proof. Note that $w_1=Aw_0-Bw_{-1}=A-\ve B$ and $(Bw_0)^2=(Aw_0-w_1)^2$.
Applying Theorem \ref{Th1.1}, we immediately get
the desired identity \eqref{sw}. \qed

\begin{corollary} {\rm (i)} For any integer $n\gs2$, we have
\begin{equation}\label{qt}\det[q^{|j-k|}+t]_{1\ls j,k\ls n}
=(1-q)^{n-1}(1+q)^{n-2}((n(1-q)+2q)t+q+1).
\end{equation}

{\rm (ii)} For any positive integer $n$, we have
\begin{equation}\label{qjk}\det[q^{|j-k|}-q^j-q^k+1]_{1\ls j,k\ls n}=(1-q^2)^n+n(1+q)^{n-1}(1-q)^{n+1}.
\end{equation}
\end{corollary}
\Proof. (i) Let $w_n=q^n+t$ for $n\in\N$. Then $w_0=t+1$, $w_1=q+t$, and
$$w_{n+1}=(q+1)w_n-qw_{n-1}\ \t{for all}\ n=1,2,3,\ldots.$$
Thus, applying Theorem \ref{Th1.1} we get the desired identity \eqref{qt}.

{\rm (ii)} By \cite[Lemma 2.1]{S19},
\begin{equation}\label{M}
\det[q^{|j-k}+t]_{1\ls j,k\ls n+1}=\det[q^{|j-k|}]_{1\ls j,k\ls n+1}+t\det(M),
\end{equation}
where $M=[m_{j,k}]_{2\ls j,k\ls n+1}$ with $m_{j,k}=q^{|j-k|}-q^{|j-1|}-q^{|1-k|}+q^{|1-1|}$.
In view of Theorem \ref{Th1.1}(i) and the identity \eqref{M},
\begin{align*}\det(M)&=(1-q)^n(1+q)^{n-1}((n+1)(1-q)+2q)
\\&=(1-q)^n(1+q)^{n-1}(n(1-q)+1+q)
\\&=n(1-q)^{n+1}(1+q)^{n-1}+(1-q)^n(1+q)^n.
\end{align*}
Note that $M=[m_{j+1,k+1}]_{1\ls j,k\ls n}$ and
$$m_{j+1,k+1}=q^{|j-k|}-q^j-q^k+1\ \ \t{for all}\ j,k=1,\ldots,n.$$
So we have the identity \eqref{qjk}. This ends our proof.
\qed

In contrast with Theorem \ref{Th1.1}, we also have the following (relatively easier) result.

\begin{theorem}\label{Th1.2} Let $A$ and $B$ be elements of a commutative ring $R$ with identity, and let $(w_n)_{n\gs0}$ be a sequence of elements of $R$
satisfying the recurrence \eqref{w}. For any integer $n>1$, we have
\begin{equation}\det[w_{|j-k+1|}]_{1\ls j,k\ls n}=(w_1^2-Aw_0w_1+Bw_0^2)((B+1)w_1-Aw_0)^{n-2}.
\end{equation}
\end{theorem}

Clearly, Theorem \ref{Th1.2} has the following consequence.

\begin{corollary} Let $A$ and $B$ be elements of a commutative ring $R$ with identity.  For any integer $n>1$, we have
\begin{equation}\det[u_{|j-k+1|}(A,B)]_{1\ls j,k\ls n}=(B+1)^{n-2}
\end{equation}
and
\begin{equation}\det[v_{|j-k+1|}(A,B)]_{1\ls j,k\ls n}=(4B-A^2)(A(B-1))^{n-2}.
\end{equation}
\end{corollary}

Now, we present our third theorem.

\begin{theorem} \label{Th1.3} For any positive integer $n$, we have
\begin{equation}\label{qx}\begin{aligned}\det[q^{|j-k|}+x\da_{jk}]_{1\ls j,k\ls n}
=\ &(x+1)u_n(1-q^2+(1+q^2)x,q^2x^2)
\\\ &-q^2x^2u_{n-1}(1-q^2+(1+q^2)x,q^2x^2),
\end{aligned}
\end{equation}
where the Kronecker symbol $\da_{jk}$ is $1$ or $0$ according as $j=k$ or not.
\end{theorem}
Taking $x=\pm1$ in \eqref{qx} and recalling the identity \eqref{xyz}, we obtain the following corollary.

\begin{corollary} \label{Cor-q+}Let $n$ be any positive integer. Then
\begin{equation}\label{q+}\det[q^{|j-k|}+\da_{jk}]_{1\ls j,k\ls n}=u_{n+1}(2,q^2)
\end{equation}
and
\begin{equation}\det[q^{|j-k|}-\da_{jk}]_{1\ls j,k\ls n}=(-1)^{n-1}q^nu_{n-1}(2q,1).
\end{equation}
\end{corollary}

Applying \eqref{q+} with $q=2$, and noting that
$$u_{n+1}(2,4)=2^nu_{n+1}(1,1)=2^n(-1)^{n}\l(\f {n+1}3\r)$$
with the aid of \eqref{xyz}, we get from Corollary \ref{Cor-q+} the following consequence.

\begin{corollary} For any positive integer $n$, we have the identity
\begin{equation}\det[2^{|j-k|}+\da_{jk}]_{1\ls j,k\ls n}=(-2)^n\l(\f{n+1}3\r).
\end{equation}
\end{corollary}

Now we state our last theorem.

\begin{theorem} \label{Th1.4} Let $A\in\C$ with $A(A^2+4)\not=0$. And let $(w_i)_{i\gs0}$ be a sequence of complex numbers
with $w_{i+1}=Aw_i+w_{i-1}$ for all $i=1,2,3,\ldots$.
For any positive even integer $n$,  $\det[x\da_{jk}-w_{j+k}]_{0\ls j,k\ls n-1}$
$($the characteristic polynomial of the matrix
$W=[w_{j+k}]_{0\ls j,k\ls n-1})$ equals
\begin{equation}\begin{aligned}x^n&-(w_1v_{n-1}(A,-1)+w_0v_{n-2}(A,-1))\f{u_n(A,-1)}{A}x^{n-1}
\\\ &+(w_0^2+Aw_0w_1-w_1^2)\f{u_n(A,-1)^2}{A^2}x^{n-2}.
\end{aligned}
\end{equation}
\end{theorem}

Taking $w_0=0$ and $w_1=1$ in Theorem \ref{Th1.4}, we get the following corollary.

\begin{corollary} \label{un-1}Let $A\in\C$ with $A(A^2+4)\not=0$. 
For any positive even integer $n$,  we have
\begin{equation}\begin{aligned}&\ A^2\det[x\da_{jk}-u_{j+k}(A,-1)]_{0\ls j,k\ls n-1}
\\=&\ A^2x^n-Au_n(A,-1)v_{n-1}(A,-1)x^{n-1}-u_n(A,-1)^2x^{n-2}.
\end{aligned}
\end{equation}
\end{corollary}

Applying Theorem \ref{Th1.4} with $w_i=u_{i+2}(A,-1)$ for all $i\in\N$, 
and noting that 
\begin{align*}v_{n+1}(A,-1)=\ &Av_n(A,-1)+v_{n-1}(A,-1)
\\=\ &A(Av_{n-1}(A,1)+v_{n-2}(A,-1))+v_{n-1}(A,-1)
\\=\ &(A^2+1)v_{n-1}(A,-1)+Av_{n-2}(A,-1)
\end{align*}
for all $n=2,3,\ldots$, we obtain the following corollary.

\begin{corollary}\label{u1n} Let $A\in\C$ with $A(A^2+4)\not=0$.
For any positive even integer $n$,  we have
\begin{equation}\begin{aligned}&\ A^2\det[x\da_{jk}-u_{j+k}(A,-1)]_{1\ls j,k\ls n}
\\=&\ A^2x^n-Au_n(A,-1)v_{n+1}(A,-1)x^{n-1}-u_n(A,-1)^2x^{n-2}.
\end{aligned}
\end{equation}
\end{corollary}

Taking $w_0=2$ and $w_1=A$ in Theorem \ref{Th1.4}, and noting that
$$Av_m(A,-1)+2v_{m-1}(A,-1)=(A^2+4)u_m(A,-1)\ \t{for all}\ m=1,2,3,\ldots$$
(which can be easily proved by induction),
we obtain the following result.

\begin{corollary} \label{vn-1} Let $A\in\C$ with $A(A^2+4)\not=0$.
For any positive even integer $n$,  we have
\begin{equation}\begin{aligned}&\ A^2\det[x\da_{jk}-v_{j+k}(A,-1)]_{0\ls j,k\ls n-1}
\\=&\ A^2x^n-A(A^2+4)u_n(A,-1)u_{n-1}(A,-1)x^{n-1}+(A^2+4)u_n(A,-1)^2x^{n-2}.
\end{aligned}
\end{equation}
\end{corollary}

Applying Theorem \ref{Th1.4} with $w_i=v_{i+2}(A,-1)$ for all $i\in\N$,
and noting that
\begin{align*}A(A^2+3)v_{n-1}(A,-1)+(A^2+2)v_{n-2}(A,-1)=(A^2+4)u_{n+1}(A,-1)
\end{align*}
for all $n=2,3,\ldots$ (which can be easily proved by induction), we get the following corollary.

\begin{corollary} \label{v1n} Let $A\in\C$ with $A(A^2+4)\not=0$.
For any positive even integer $n$,  we have
\begin{equation}\begin{aligned}&\ A^2\det[x\da_{jk}-v_{j+k}(A,-1)]_{1\ls j,k\ls n}
\\=&\ A^2x^n-A(A^2+4)u_n(A,-1)u_{n+1}(A,-1)x^{n-1}+(A^2+4)u_n(A,-1)^2x^{n-2}.
\end{aligned}
\end{equation}
\end{corollary}

The Lucas numbers are those $L_n=v_n(1,-1)\ (n\in\N)$.
Taking $x=-1$ in Corollaries \ref{u1n} and \ref{v1n}, we obtain the following consequence.

\begin{corollary} Let $n$ be a positive even number. For any $A\in\C$ with $A(A^2+4)\not=0$, we have
\begin{equation}A^2\det[u_{j+k}(A,-1)+\da_{jk}]_{1\ls j,k\ls n}=(A-1)(A+u_n(A,-1)^2)+Au_{n+1}(A,-1)^2
\end{equation}
and
\begin{equation}A^2\det[v_{j+k}(A,-1)+\da_{jk}]=v_{n+1}(A,-1)^2.\end{equation}
 In particular,
\begin{equation}\det[F_{j+k}+\da_{jk}]_{1\ls j,k\ls n}=F_{n+1}^2
\end{equation}
and
\begin{equation}\det[L_{j+k}+\da_{jk}]_{1\ls j,k\ls n}=L_{n+1}^2.
\end{equation}
\end{corollary}

Similarly, taking $x=-1$ and $A=1$ in Corollaries \ref{un-1} and \ref{vn-1}, we find that for any positive even integer $n$ we have
\begin{equation}\det[F_{j+k}+\da_{jk}]_{0\ls j,k\ls n-1}=F_{n-1}^2
\end{equation}
and
\begin{equation}\det[L_{j+k}+\da_{jk}]_{0\ls j,k\ls n-1}=L_nL_{n+1}-1=L_{2n+1}.
\end{equation}

We are going to prove Theorems \ref{Th1.1}-\ref{Th1.2} and Theorems \ref{Th1.3}-\ref{Th1.4} in Sections 2 and 3, respectively.
We will propose some conjectures in Section 4.

\section{Proofs of Theorems \ref{Th1.1}-\ref{Th1.2}}
 \setcounter{equation}{0}
 \setcounter{conjecture}{0}
 \setcounter{theorem}{0}
 \setcounter{proposition}{0}

\medskip
 \noindent{\it Proof of Theorem \ref{Th1.1}}.
 Let $M_n$ denote the matrix $[w_{|j-k|}]_{1\ls j,k\ls n}$. Clearly,
 $$\det(M_1)=w_0=w_0u_1(A',B')+((Bw_0)^2-(w_1-Aw_0)^2)u_0(A',B')$$
 and
 \begin{align*}\det(M_2)&=\begin{vmatrix}
w_0&w_1\\w_1&w_0
\end{vmatrix}=w_0^2-w_1^2
\\&=w_0u_2(A',B')+((Bw_0)^2-(w_1-Aw_0)^2)u_1(A',B').
\end{align*}
So \eqref{w'} holds for $n=1,2$.

Now suppose $n\gs3$, and assume that
\begin{equation}\label{k}\det(M_{k})=w_0u_{k}(A',B')+((Bw_0)^2-(w_1-Aw_0)^2)u_{k-1}(A',B')
\end{equation}
for each $k=1,\ldots,n-1$.
Observe that
 \begin{align*}&Bw_{|(n-2)-k|}-Aw_{|(n-1)-k|}+w_{|n-k|}
 \\=\ &\begin{cases}Bw_{n-2-k}-Aw_{n-1-k}+w_{n-k}=0&\t{if}\ 1\ls k<n-1,
 \\Bw_1-Aw_0+w_1=(B+1)w_1-Aw_0&\t{if}\ k=n-1,
 \\Bw_2-Aw_1+w_0=B(Aw_1-Bw_0)-Aw_1+w_0&\t{if}\ k=n.
 \end{cases}
 \end{align*}
Thus, adding the $(n-2)$-th row times $B$ and the $(n-1)$-th row times $-A$ to the last row
 of $M_n$, we find that $\det(M_n)=\det(M_n')$, where
 $$M_n':=\begin{bmatrix}
w_0&w_1&w_2&\cdots&w_{n-3}&w_{n-2}&w_{n-1}\\
w_1&w_0&w_1&\cdots&w_{n-4}&w_{n-3}&w_{n-2}\\
\vdots&\vdots&\vdots&\ddots&\vdots&\vdots&\vdots\\
w_{n-4}&w_{n-5}&w_{n-6}&\cdots&w_1&w_2&w_3\\
w_{n-3}&w_{n-4}&w_{n-5}&\cdots&w_0&w_1&w_2\\
w_{n-2}&w_{n-3}&w_{n-4}&\cdots&w_1&w_0&w_1\\
0&0&0&\cdots&0&C&(B-1)D
\end{bmatrix}$$
with $C =(B+1)w_1-Aw_0$ and $D=Aw_1-(B+1)w_0$.
Adding the $(n-2)$-th column times $B$ and the $(n-1)$-th column times $-A$ to the last column
 of $M_n'$, we see that $\det(M_n')=\det(M_n'')$, where
 $$M_n'':=\begin{bmatrix}
w_0&w_1&w_2&\cdots&w_{n-3}&w_{n-2}&0\\
w_1&w_0&w_1&\cdots&w_{n-4}&w_{n-3}&0\\
\vdots&\vdots&\vdots&\ddots&\vdots&\vdots&\vdots\\
w_{n-4}&w_{n-5}&w_{n-6}&\cdots&w_1&w_2&0\\
w_{n-3}&w_{n-4}&w_{n-5}&\cdots&w_0&w_1&0\\
w_{n-2}&w_{n-3}&w_{n-4}&\cdots&w_1&w_0&C\\
0&0&0&\cdots&0&C&(B-1)D-AC
\end{bmatrix}.$$
Expanding $\det(M_n'')$ via its last row, we get
\begin{align*}\det(M_n'')=\ &((B-1)D-AC)\begin{vmatrix}
w_0&w_1&w_2&\cdots&w_{n-3}&w_{n-2}\\
w_1&w_0&w_1&\cdots&w_{n-4}&w_{n-3}\\
\vdots&\vdots&\vdots&\ddots&\vdots&\vdots\\
w_{n-4}&w_{n-5}&w_{n-6}&\cdots&w_1&w_2\\
w_{n-3}&w_{n-4}&w_{n-5}&\cdots&w_0&w_1\\
w_{n-2}&w_{n-3}&w_{n-4}&\cdots&w_1&w_0
\end{vmatrix}
\\&-C\begin{vmatrix}
w_0&w_1&w_2&\cdots&w_{n-3}&0\\
w_1&w_0&w_1&\cdots&w_{n-4}&0\\
\vdots&\vdots&\vdots&\ddots&\vdots&\vdots\\
w_{n-4}&w_{n-5}&w_{n-6}&\cdots&w_1&0\\
w_{n-3}&w_{n-4}&w_{n-5}&\cdots&w_0&0\\
w_{n-2}&w_{n-3}&w_{n-4}&\cdots&w_1&C\\
\end{vmatrix}.\end{align*}
Therefore
$$\det(M_n)=\det(M_n')=\det(M_n'')=((B-1)D-AC)\det(M_{n-1})-C^2\det(M_{n-2}).$$
Note that $C^2=B'$ and
$$(B-1)D-AC=(B-1)(Aw_1-(B+1)w_0)-A((B+1)w_1-Aw_0)=A'.$$
Thus, with the aid of \eqref{k} for $k=n-1,n-2$, we have
\begin{align*}\det(M_n)=\ &A'\det(M_{n-1})-B'\det(M_{n-2})
\\=\ &A'(w_0u_{n-1}(A',B')+((Bw_0)^2-(w_1-Aw_0)^2)u_{n-2}(A',B'))
\\\ &-B'(w_0u_{n-2}(A',B')+((Bw_0)^2-(w_1-Aw_0)^2)u_{n-3}(A',B'))
\\=\ &w_0u_n(A',B')+((Bw_0)^2-(w_1-Aw_0)^2)u_{n-1}(A',B').
\end{align*}

In view of the above, by induction the identity \eqref{w'} holds for any $n\in\Z^+$. \qed

 \medskip
 \noindent{\it Proof of Theorem \ref{Th1.2}}.
 Let $W_n=\det[w_{|j-k+1|}]_{1\ls j,k\ls n}$. Clearly,
 $$W_2=\begin{vmatrix}
w_1&w_0\\w_2&w_1
\end{vmatrix}=w_1^2-w_0(Aw_1-Bw_0)=w_1^2-Aw_0w_1+Bw_0^2.$$

Now, assume $n\gs3$.
 Observe that
 \begin{align*}&Bw_{|(n-2)-k+1|}-Aw_{|(n-1)-k+1|}+w_{|n-k+1|}
 \\=\ &\begin{cases}Bw_{n-k-1}-Aw_{n-k}+w_{n-k+1}=0&\t{if}\ 1\ls k<n,
 \\Bw_1-Aw_0+w_1=(B+1)w_1-Aw_0&\t{if}\ k=n.
 \end{cases}
 \end{align*}
 Thus, adding the $(n-2)$-th row times $B$ and the $(n-1)$-th row times $-A$ to the last row
 of the determinant $W_n$, we find that the last row turns to be
 $$\underbrace{0\ \cdots\ 0}_{n-1}\ (B+1)w_1-Aw_0. $$
 Therefore
 $$W_n=((B+1)w_1-Aw_0)W_{n-1}.$$

 In view of the above, by induction we have
$$W_n=(w_1^2-Aw_0w_1+Bw_0^2)((B+1)w_1-Aw_0)^{n-2}$$
for all $n=2,3,\ldots$. This ends our proof. \qed

\section{Proofs of Theorems \ref{Th1.3}-\ref{Th1.4}}
 \setcounter{equation}{0}
 \setcounter{conjecture}{0}
 \setcounter{theorem}{0}
 \setcounter{proposition}{0}

 \medskip
 \noindent{\it Proof of Theorem \ref{Th1.3}}. Let $a_{jk}=q^{|j-k|}+x\da_{jk}$ for all $j,k=1,\ldots,n$, and let
 $Q_n$ denote the matrix $[a_{jk}]_{1\ls j,k\ls n}$. Clearly,
 $\det(Q_1)=q^0+x=x+1$
 and $$\det(Q_2)=\begin{vmatrix}
q^0+x&q\\q&q^0+x\end{vmatrix}=(x+1)^2-q^2=(x+1)(1-q^2+(1+q^2)x)-q^2x^2.$$
Thus, $\eqref{qx}$ holds for $n\in\{1,2\}$.

Now, we let $n\gs 3$ and assume the equality
\begin{equation}\label{m}\begin{aligned}\det[q^{|j-k|}+x\da_{jk}]_{1\ls j,k\ls m}
=\ &(x+1)u_m(1-q^2+(1+q^2)x,q^2x^2)
\\\ &-q^2x^2u_{m-1}(1-q^2+(1+q^2)x,q^2x^2)
\end{aligned}
\end{equation}
for any positive integer $m<n$.

 Observe that
 \begin{align*}-qa_{n-1,k}+a_{nk}&=-q(q^{|n-1-k|}+x\da_{n-1,k})+q^{|n-k|}+x\da_{nk}
 \\&=\begin{cases}0&\t{if} 1\ls k\ls n-2,\\-q(1+x)+q=-qx&\t{if}\ k=n-1,
 \\-q^2+1+x&\t{if}\ k=n.\end{cases}
 \end{align*}
 Thus, via adding the $(n-1)$-th row times $-q$ to the $n$-th row of $Q_n$, we see that
 $\det(Q_n)=\det(Q_n')$, where
 $$Q_n':=\begin{bmatrix}
a_{11}&a_{12}&a_{13}&\cdots&a_{1,n-2}&a_{1,n-1}&a_{1,n}\\
a_{21}&a_{22}&a_{23}&\cdots&a_{2,n-2}&a_{2,n-1}&a_{2,n}\\
\vdots&\vdots&\vdots&\ddots&\vdots&\vdots&\vdots\\
a_{n-3,1}&a_{n-3,2}&a_{n-3,3}&\cdots&a_{n-3,n-2}&a_{n-3,n-1}&a_{n-3,n}\\
a_{n-2,1}&a_{n-2,2}&a_{n-2,3}&\cdots&a_{n-2,n-2}&a_{n-2,n-1}&a_{n-2,n}\\
a_{n-1,1}&a_{n-1,2}&a_{n-1,3}&\cdots&a_{n-1,n-2}&a_{n-1,n-1}&a_{n-1,n}\\
0&0&0&\cdots&0&-qx&x+1-q^2
\end{bmatrix}.$$
Via adding the $(n-1)$-th column times $-q$ to the $n$-th column of $Q_n'$, we get that
 $\det(Q_n')=\det(Q_n'')$, where
 $$Q_n'':=\begin{bmatrix}
a_{11}&a_{12}&a_{13}&\cdots&a_{1,n-2}&a_{1,n-1}&0\\
a_{21}&a_{22}&a_{23}&\cdots&a_{2,n-2}&a_{2,n-1}&0\\
\vdots&\vdots&\vdots&\ddots&\vdots&\vdots&\vdots\\
a_{n-3,1}&a_{n-3,2}&a_{n-3,3}&\cdots&a_{n-3,n-2}&a_{n-3,n-1}&0\\
a_{n-2,1}&a_{n-2,2}&a_{n-2,3}&\cdots&a_{n-2,n-2}&a_{n-2,n-1}&0\\
a_{n-1,1}&a_{n-1,2}&a_{n-1,3}&\cdots&a_{n-1,n-2}&a_{n-1,n-1}&-qx\\
0&0&0&\cdots&0&-qx&f(q,x)
\end{bmatrix},$$
where
$$f(q,x):=1-q^2+(1+q^2)x.$$
Expanding $\det(Q_n'')$ via its last row, we see that
\begin{align*}\det(Q_n'')=&\ f(q,x)\det(Q_{n-1})
\\&\  -(-qx)\begin{vmatrix}
a_{11}&a_{12}&a_{13}&\cdots&a_{1,n-2}&0\\
a_{21}&a_{22}&a_{23}&\cdots&a_{2,n-2}&0\\
\vdots&\vdots&\vdots&\ddots&\vdots&\vdots\\
a_{n-3,1}&a_{n-3,2}&a_{n-3,3}&\cdots&a_{n-3,n-2}&0\\
a_{n-2,1}&a_{n-2,2}&a_{n-2,3}&\cdots&a_{n-2,n-2}&0\\
a_{n-1,1}&a_{n-1,2}&a_{n-1,3}&\cdots&a_{n-1,n-2}&-qx
\end{vmatrix}.
\end{align*}
Therefore,
$$\det(Q_n)=\det(Q_n')=\det(Q_n'')=f(q,x)\det(Q_{n-1})-(-qx)^2\det(Q_{n-2}),$$
Combining this with \eqref{m} for $m=n-1,n-2$, we find that
\begin{align*}\det(Q_n)=\ &f(q,x)\l((x+1)u_{n-1}(f(q,x),q^2x^2)-q^2x^2u_{n-2}(f(q,x),q^2x^2)\r)
\\\ &-q^2x^2\l((x+1)u_{n-2}(f(q,x),q^2x^2)-q^2x^2u_{n-3}(f(q,x),q^2x^2)\r)
\\=\ &(x+1)u_n(f(x,q),q^2x^2)-q^2x^2u_{n-1}(f(q,x),q^2x^2).
\end{align*}

In view of the above, we have proved the desired result by induction $n$. \qed

\begin{lemma}\label{Lem-w} Let $A,B\in\C$, and let
$$w_0,w_1\in\C,\ \t{and}\ w_{i+1}=Aw_i-Bw_{i-1}\ \t{for all}\ i=1,2,3,\ldots.$$
Then, for any $j\in\N$ and $k\in\Z^+$, we have
\begin{equation}w_{j+k}=w_{j+1}u_k(A,B)-Bw_ju_{k-1}(A,B).
\end{equation}
\end{lemma}
\Proof. This can be easily proved by induction on $k$. \qed

\medskip
\noindent{\it Proof of Theorem \ref{Th1.4}}. Let $\al$ and $\beta$ be the two distinct roots of the quadratic equation
$x^2-Ax+B=0$ with $B=-1$. Then $\al+\beta=A$ and $\al\beta=B=-1$.
It is well known that there are $a,b\in\C$ such that $w_m=a\al^m+b\beta^m$ for all $m\in\Z$.
As $a+b=w_0$ and $a\al+b\beta=w_1$, we find that
\begin{equation}\label{a,b}a=\f{w_1-\beta w_0}{\al-\beta}\ \ \t{and}\ \ \ b=\f{\al w_0-w_1}{\al-\beta}.
\end{equation}
It follows that
\begin{equation}\label{ab}ab=\f{(\al+\beta)w_0w_1-Bw_0^2-w_1^2}{(\al-\beta)^2}=\f{w_0^2+Aw_0w_1-w_1^2}{(\al-\beta)^2}.
\end{equation}
Observe that $\al\not=\pm1$ since $\al\beta=B=-1$ and $\al+\beta=A\not=0$.
For any $j\in\N$, we clearly have
\begin{align*}\sum_{k=0}^{n-1}w_{j+k}\al^k&=\sum_{k=0}^{n-1}(a\al^{j+k}+b\beta^{j+k})\al^k
\\&=a\al^j\sum_{k=0}^{n-1}\al^{2k}+b\beta^j\sum_{k=0}^{n-1}B^k=a\al^j\f{\al^{2n}-1}{\al^2-1}+b\beta^j\f{B^n-1}{B-1}.
\end{align*}
Note that $B^n=1$ as $B=-1$ and $2\mid n$. Thus
$$\sum_{k=0}^{n-1}w_{j+k}\al^k=a\la^j\f{\al^{2n}-(\al\beta)^n}{\al^2+\al\beta}=\al^j\times\f{a}A\al^{n-1}(\al^n-\beta^n).$$
So, the matrix $W=[w_{j+k}]_{0\ls j,k\ls n-1}$ has an eigenvalue
$$\la_0=\f aA\al^{n-1}(\al^n-\beta^n)$$
with the eigenvector $(1,\al,\al^2,\ldots,\al^{n-1})^T$.
Similarly,
$$\la_1=\f bA\beta^{n-1}(\beta^n-\al^n)=-\f bA\beta^{n-1}(\al^n-\beta^n)$$
is an eigenvalue of $W$ with the eigenvector $(1,\beta,\beta^2,\ldots,\beta^{n-1})^T$.
If $c,d\in\C$ are not all zero, and $c\al^k+d\beta^k=0$
for all $k=0,\ldots,n-1$, then $c=-d\not=0$ and $\al=\beta$. As $\al\not=\beta$, the two vectors
$(1,\al,\al^2,\ldots,\al^{n-1})^T$ and $(1,\beta,\beta^2,\ldots,\beta^{n-1})^T$ are linearly independent over $\C$.

For each $2\ls m\ls n-1$, since
$$w_n(-Bu_{m-1})+w_{j+1}u_m+w_{j+m}\times(-1)=0$$
by Lemma \ref{Lem-w}, the vector $V_m=(v_{m0},v_{m1},\ldots,v_{m,n-1})^T$ is an eigenvector associated to the eigenvalue
$\la_{m}=0$ of the matrix $W$, where
$$v_{m0}=-Bu_{m-1}=u_{m-1},\ v_{m1}=u_m,\ v_{mm}=-1,$$
 and $v_{mk}=0$ for all $2\ls k\ls n-1$ with $k\not=m$.
 If $c_2,\ldots,c_{n-1}\in\C$ and $\sum_{m=2}^{n-1}c_mv_{mk}=0$ for all $k=0,1,\ldots,n-1$,
 then for each $2\ls k\ls n-1$ we have
 $$0=\sum_{m=2}^{n-1}c_mv_{mk}=c_kv_{kk}=c_k.$$
 So, the vectors $V_2,\ldots,V_{n-1}$ are linearly independent over $\C$.

In view of the above, $\la_0,\la_1,\la_2,\ldots,\la_{n-1}$ are all the $n$ eigenvalues of the matrix $W=[w_{j+k}]_{0\ls j,k\ls n-1}$. So we have
\begin{align*}\det[x\da_{jk}-w_{j+k}]_{0\ls j,k\ls n-1}&=\prod_{m=0}^{n-1}(x-\lambda_m)
=x^{n-2}(x-\la_0)(x-\la_1)\\&=x^{n-2}(x^2-(\la_0+\la_1)x+\la_0\la_1).
\end{align*}
In view of \eqref{a,b},
$$\la_0+\la_1=\f{u_n}A\l((w_1-\beta w_0)\al^{n-1}-(\al w_0-w_1)\beta^{n-1}\r)=\f{u_n}A\l(w_1v_{n-1}-Bw_0v_{n-2}\r).$$
Also,
$$\la_0\la_1=-\f{abB^{n-1}}{A^2}(\al^n-\beta^n)^2=-\f{B^{n-1}u_n^2}{A^2}(w_0^2+Aw_0w_1-w_1^2)$$
by \eqref{ab}.

\section{Some Conjectures}
 \setcounter{equation}{0}
 \setcounter{conjecture}{0}
 \setcounter{theorem}{0}
 \setcounter{proposition}{0}

 \begin{conjecture} For any positive integer $n$, we have
 \begin{equation}\det[2^{|j-k|}-1+\da_{j,k}]_{1\ls j,k\ls n}=2^n+(-1)^n2^{n-1}\l(2\l(\f n3\r)+n\l(\f{n+1}3\r)\r).
 \end{equation}
 \end{conjecture}

 \begin{remark} For any positive integer $n$,  \cite[Theorem 1.4]{WS} implies  that
 $$\det[2^{j+k}-1+\da_{jk}]_{1\ls j,k\ls n}=4(2^n-1)^2-(n-1)\f{4^{n+1}-1}3.$$
 \end{remark}

 \begin{conjecture} For any positive odd integer $n$, we have
 \begin{align}\det[F_{j+k}+\da_{jk}]_{0\ls j,k\ls n-1}&=F_{n+1}^2+1
 \\\det[F_{j+k}+\da_{jk}]_{1\ls j,k\ls n}&=F_{n+1}^2+1=F_nF_{n+2},
 \\\det[L_{j+k}+\da_{jk}]_{0\ls j,k\ls n-1}&=L_nL_{n+1}=L_{2n+1}-1.
 \end{align}
 \end{conjecture}

 \begin{remark} For any positive integer $n$, H. Wang and Z.-W. Sun \cite[Theorem 1.1(ii)]{WS} proved that
 $$\det[F_{|j-k|}+\da_{jk}]_{1\ls j,k\ls n}=\begin{cases}1&\t{if}\ n\eq0,\pm1\pmod6,\\0&\t{otherwise}.
 \end{cases}$$
 \end{remark}

 Based on our computation via {\tt Mathematica}, we also propose the following two conjectures.

 \begin{conjecture} For any $A\in\C$ and $n\in\Z^+$, we have
 \begin{equation}\det[v_{j+k}(A,1)+\da_{jk}]_{1\ls j,k\ls n}=u_{n+1}(A,1)^2-n^2.
 \end{equation}
 \end{conjecture}

\begin{conjecture} Let $n$ be any positive odd integer. For any $A\in\C$ we have
 \begin{equation}A^2\det[v_{j+k}(A,-1)+\da_{jk}]_{1\ls j,k\ls n}= v_{n+1}(A,-1)^2-A^2-4.
 \end{equation}
 In particular,
 \begin{equation}\det[L_{j+k}+\da_{jk}]_{1\ls j,k\ls n}= L_{n+1}^2-5.
\end{equation}
 \end{conjecture}
 \begin{remark} For any $A,B\in\C$ and $n\in\N$, it is well known that
 $$v_n(A,B)^2-(A^2-4B)u_n(A,B)^2=4B^n$$
 which can be easily proved.
 \end{remark}

\end{document}